\documentclass[review]{elsarticle}

\usepackage{lineno,hyperref}
\modulolinenumbers[5]
\usepackage{amsmath}
\usepackage{amssymb}
 \usepackage{paralist,float,lscape}
 \usepackage{movie15}
 
\newcommand{\R}{{\mathbb  R}}
\newtheorem{theorem}{Theorem}[section]
\newtheorem{remark}{Remark}[section]
\newtheorem{prop}{Proposition}[section]
\newtheorem{definition}{Definition}[section]
\journal{ }

\begin{document}

\begin{frontmatter}

\title{Perturbation of periodic orbits of Rabinovich system}

\author{G\^\i rban Anania\fnref{myfootnote}}
\address{B-dul. Vasile P\^{a}rvan, No. 4, 300223 - Timi\c soara, Rom\^ania}

\author[anania.girban@upt.ro]{Politehnica University of Timi\c soara}

\cortext[mycorrespondingauthor]{Corresponding author}
\ead{anania.girban@upt.ro}

\begin{abstract}
In this paper we stabilize asymptotically the periodic orbits of the Rabinovich system.
\end{abstract}

\begin{keyword}
\texttt{Hamiltonian dynamics, Rabinovich system,  periodic orbits}
\MSC[2010] 00-01\sep  99-00
\end{keyword}

\end{frontmatter}


\section{Introduction}\label{section:one}

The Rabinovich system (see e.g. \cite{rabinovich}, \cite{rabinovich_1}) has been intensely studied last time from many points of view 
(see e.g.  \cite{xie}, \cite{zhang}, \cite{llibre_1}, \cite{razvan} and many others). In  \cite{razvan} we studied this system from the Poisson geometry and the dynamics point of view. 
We studied in particularly the existence of periodic solutions by using the Lyapunov center theorem on the symplectic leaves of the Poisson configuration manifold.

The purpose of this article is to stabilize asymptotically the periodic orbits of the Rabinovich system using a method described in \cite{razvan1} and to find a perturbation which stabilize 
globally asymptotically an arbitrary periodic orbit of a generic three-dimensional Hamiltonian system and in particular the periodic orbits of the Rabinovich system.
In the paper \cite{razvan} the Rabinovich system is modeled as a Hamilton-Poisson dynamical system on a Poisson manifold with one Casimir and consequently, the dynamics is given by common level 
set of  the Casimir and the Hamiltonian of the Poisson manifold.

We consider an arbitrary periodic orbit of the Rabinovich system around a Lyapunov stable equilibrium point. For this orbit we don't know any parameterization 
or if it is asymptotically stable or unstable. From the Poisson structure, we know that the periodic orbit of the system lies on a common level set of the Casimir and the Hamiltonian. 
Now we fix this common level set of the Casimir and the Hamiltonian and perturb the system, but the periodic orbit remains a periodic orbit for the perturbed system too. The perturbed families are 
parameterized by an arbitrary strictly positive smooth real function.
There are two cases, one when the Hamiltonian remains dynamically invariant for the perturbed system and   one when the Casimir remains  dynamically invariant for the perturbed system. 
We will consider and a third case when the periodic orbit remains a periodic orbit for the perturbed system too, but the Hamiltonian and the Casimir do not remain dynamically invariant for the 
perturbed system.
In the first two cases we have two possibilities: first when the periodic orbit of the perturbed system can be the orbitally phase asymptotically stable and second when it is unstable with respect 
to  perturbations along the level set of the Hamiltonian (respectively the Casimir), which contains the periodic orbit.

But, in the first two cases, we can't obtain the global asymptotic stability because the dynamics is located on the common level surfaces Hamiltonian=constant or Casimir=constant and at least one 
of them is dynamically invariant in this cases. If we do not keep both the Hamiltonian and the Casimir dynamically invariant, only the periodic orbit of the studied system remains periodic orbit for 
the perturbed system, we can  stabilize asymptotically this periodic orbit in an entire tubular neighborhood around it.
So, we consider a third case for a generic three-dimensional Hamiltonian system for global asymptotic stabilization of a periodic orbit of the system when the Hamiltonian and the Casimir do not remain 
dynamically invariant for the perturbed system and consequently the dynamics is located in an entire tubular neighborhood around the periodic orbit of the system without the constrains to remain 
on the surfaces Hamiltonian=constant or Casimir=constant from the first two cases.
  In this third case, the periodic orbit of the given system remains periodic orbit for the perturbed system and the perturbed families are parameterized by two arbitrary strictly positive smooth 
  real functions.

\section{Preliminaries}
The Rabinovich system which we propose to study, namely
\begin{equation}\label{sys}
\left\{ \begin{array}{l}
 \dot x = yz + \beta y \\
 \dot y =-xz + \beta x \\
 \dot z = xy \\
 \end{array} \right.
\end{equation}
where $\beta\in\R$ is a parameter, is a particular case of the Rabinovich system introduced in \cite{rabinovich}.

In the paper  \cite{razvan} we studied the system \eqref{sys} from Poisson geometry and the dynamics point of view. 
We recall from \cite{razvan} some results that we need to stabilize asymptotically the periodic orbits of the Rabinovich system. This results are: the Poisson structure (Proposition \ref{hc}), 
and  the Casimir of the configuration (Proposition \ref{hc}), because the system has to be a Hamiltonian system and in the first two cases we keep invariant first the 
Hamiltonian and then the Casimir, the equilibrium states and their stability (Proposition \ref{stab}) and the existence of periodic orbits (Proposition \ref{per_orb1} and Proposition \ref{per_orb2}) because we want to perturb 
one (or a pair) arbitrary periodic orbit from a family of periodic orbits wich appear around the Lyapunov stable equilibrium states of the Rabinovich system.

\begin{prop}[\cite{razvan}]\label{hc}

The center of the Poisson algebra $C^\infty(\R^3,\R)$ is generated by the Casimir invariant $C\in C^\infty(\R^3,\R)$, $C(x,y,z)=\dfrac{1}{2}(-x^2+y^2)+z^2$ and 
the dynamics \eqref{sys} has the following Hamilton-Poisson realization:
\begin{equation}
(\R^3,\Pi^-_{LP},H_\beta)
\end{equation}
where
\begin{equation}\label{pi}
\Pi^-_{LP}(x,y,z)=\left[ {\begin{array}{*{20}r}
   0 & {2z} & {-y}  \\
   { - 2z} & 0 & {-x }  \\
   {y } & {x} & 0  \\
\end{array}} \right]
\end{equation}
is the minus Lie-Poisson structure on $(\mathfrak o(Q))^* \cong\R^3$, and the Hamiltonian $H_\beta\in C^\infty(\R^3,\R)$ is given by $H_\beta(x,y,z)=\dfrac{1}{4}(x^2+y^2)-\beta z$.

\end{prop}

Now we can write the dynamics \eqref{sys} in the form
\begin{equation}\label{sys1}
\dfrac{du}{dt}= \nu (u) (\nabla H(u) \times \nabla C(u)) ,\ u=(x,y,z) \in\R^3
\end{equation}
with $ \nu (x,y,z)=1$.

\begin{prop}[\cite{razvan}]\label{stab}
The equilibrium states of the Rabinovich system are given as the union of the following three families:
\begin{equation}\label{eq}
\begin{split}
{\mathcal E}_{1,\beta}&=\{(M,0,\beta):M\in\R\},\\
{\mathcal E}_{2,\beta}&=\{(0,M,-\beta):M\in\R\},\\
{\mathcal E}_{3,\beta}&=\{(0,0,M):M\in\R\}
\end{split}
\end{equation}
 and the stability of the equilibrium states from the each family is the following:
\begin{enumerate}
 \item All the equilibrium states from the family ${\mathcal E}_{1,\beta}$ are nonlinearly stable except for the equilibrium $(0,0,\beta)$, which is unstable for $\beta\neq 0$.

\item All the equilibrium states from the family ${\mathcal E}_{2,\beta}$ are unstable for $\beta\neq 0$. For $\beta=0$ all the equilibrium states from the family ${\mathcal E}_{2,0}$ are 
unstable except for the origin which is nonlinearly stable.

\item Let $\beta\neq 0$, and $e_M=(0,0,M)\in {\mathcal E}_{3,\beta}$ be an arbitrary equilibrium state. The equilibrium $e_M\in {\mathcal E}_{3,\beta}$ is nonlinearly stable for $|M|>|\beta|$ and 
unstable for $|M|\leq|\beta|$. For $\beta=0$ all the equilibrium states from the family ${\mathcal E}_{3,0}$ are nonlinearly stable.
\end{enumerate}
\end{prop}

Around the Lyapunov stable equilibrium points  of the above families, we poved the existence of the periodic solutions of the Rabinovich system \eqref{sys} as it shows in the below propositions:
\begin{prop}[\cite{razvan}]\label{per_orb1}
Let $(0,0,M)\in{\mathcal E}_{3,\beta}$ be such that $|M|>|\beta|$. Then there e\-xists $\varepsilon_0>0$ and a one-parameter family 
$\left( {\gamma _\varepsilon ^{(0,0,M)} } \right)_{0 < \varepsilon  \le \varepsilon _0}$ of periodic solutions of the Rabinovich system \eqref{sys}, that approaches 
$(0,0,M)$ as $\varepsilon\to0$ with periods
$T_\varepsilon ^{(0,0,M)} \mathop{\longrightarrow}\limits^{\varepsilon\to 0}\dfrac{{2\pi }}{{\sqrt{M^2-\beta^2}}}$. The union
$$\{ (0,0,M)\}  \cup \bigcup\limits_{0 < \varepsilon\le \varepsilon _0 } {\gamma _\varepsilon ^{(0,0,M)}}$$
forms a smooth two dimensional manifold with boundary $\gamma _{\varepsilon_0}^{(0,0,M)}$, manifold that is diffeomorphic to the closed disk in $\R^2$.
\end{prop}

\begin{prop}[\cite{razvan}]\label{per_orb2}
Let $M\neq \pm\beta\sqrt{2}$ and $M\neq 0$. Then for each equilibrium state $(M,0,\beta)\in{\mathcal E}_{1,\beta}$, there exists $\delta_0 >0$ and a one-parameter family $\left({\gamma _\delta ^{(M,0,\beta)} } \right)_{0 < \delta  \le \delta _0}$ of periodic solutions of the Rabinovich system \eqref{sys}, that approaches the equilibrium $(M,0,\beta)$ as $\delta\to0$, and have the periods
$T_\delta ^{(M,0,\beta) } \mathop{\longrightarrow}\limits^{\delta\to 0}\dfrac{{2\pi }}{{|M|}}$. The union
$$\{(M,0,\beta)\}\cup \bigcup\limits_{0 < \delta\le \delta _0} {\gamma _\delta ^{(M,0,\beta)}},$$
form a two dimensional manifold with boundary $\gamma _{\delta_0}^{(M,0,\beta)}$, manifold that is diffeomorphic with the closed disk in $\R^2$.
\end{prop}

We recall now the definition concerning the stability of the periodic orbits of a dynamical system. 

Let consider  $\dot x= X(x)$ a dynamical system generated by a smooth vector field $X\in X(U)$, defined on an open subset $U \subseteq \R^n$ and
 $\Gamma=\{\gamma (t) \subset    U : 0 \leq t \leq T\}$ is a $T$-periodic orbit of $\dot x= X(x)$.
 
\begin{definition} [\cite{razvan2}]
\begin{enumerate}
 \item  The periodic orbit $\Gamma$  is called {\bf orbitally stable} if, given $\varepsilon>0$ there exists a $\delta>0$ such that 
$dist(x(t; x_0 ), \Gamma ) < \varepsilon$  for all $t > 0$ and for all $x_0 \in U$
such that $dist(x_0 , \Gamma ) < \delta$ .
\item The periodic orbit $\Gamma$ is called {\bf unstable} if it is not orbitally stable.
\item The periodic orbit $\Gamma$ is called {\bf orbitally asymptotically stable} if it is orbitally
stable and (by choosing $\delta$  smaller if necessary), $dist(x(t; x_0 ), \Gamma ) \to  0$ as $t \to\infty$.
\item The periodic orbit $\Gamma$ is called {\bf orbitally phase asymptotically stable}, if it is
orbitally asymptotically stable and there is a $\delta>0$ such that for each $x_0\in U$ with
$ dist(x_0 , \Gamma ) < \delta$, there exists $\theta_0 = \theta_0 (x_0 )$ such that
$$\lim_{t \to\infty}\|x(t; x_0 ) - \gamma (t + \theta_0 )\| = 0.$$
\end{enumerate}
\end{definition}
 For more details regarding the stability analysis of periodic orbits see e.g., \cite{razvan2}.

We recall now the main theorem from \cite{razvan1} which provides two classes of perturbations, of which one keeps dynamically invariant the Hamiltonian $H$, and respectively one keeps 
dynamically invariant the Casimir function $C$ and the periodic orbit of the initial system will remain a periodic orbit for the perturbed system too.

\begin{theorem}[\cite{razvan1}]\label{thm}
 Let
\begin{equation}\label{eq}
 \dfrac{du}{dt}= \nu (u) (\nabla H(u) \times \nabla C(u)) ,\ u \in   U
 \end{equation}
be a three-dimensional dynamical system defined on an open subset $U \subseteq\R^3$ , with
$H, C, \nu  \in   C^\infty   (U, \R)$ given smooth real functions, such that $H$ and $C$ are functionally
independent on an open subset $V \subseteq U$.

Suppose there exists $\Gamma  \subset V$ a periodic orbit of \eqref{eq}. If $\Gamma  \subseteq  (H, C)^{-1} (\{(h, c)\})$,
where $(h, c) \in   \R^2$ is a regular value for the map $(H, C) : U \to  \R^2$, then the following
conclusions hold true.

\begin{enumerate}
 \item If $c$ is a regular value of the map $C : U \to  \R$, then for every smooth function
$\alpha  \in   C^\infty   (V, (0, \infty))$:
\begin{enumerate}
 \item $\Gamma$ , as a periodic orbit of the perturbed dynamical system
$$\dfrac{du}{dt}= \nu (u) (\nabla H(u) \times \nabla C(u)) - \alpha (u)(H(u) - h) [\nabla C(u) \times (\nabla H(u) \times \nabla C(u))] ,$$
$u \in   V $, is orbitally phase asymptotically stable, with respect to perturbations
in $V$ , along the invariant manifold $C^{-1} (\{c\})$.
\item $\Gamma$ , as a periodic orbit of the perturbed dynamical system
$$\dfrac{du}{dt}= \nu (u) (\nabla H(u) \times \nabla C(u)) + \alpha (u)(H(u) - h) [\nabla C(u) \times (\nabla H(u) \times \nabla C(u))] ,$$
$u \in   V $, is unstable.
\end{enumerate}

\item If $h$ is a regular value of the map $H : U \to\R$, then for every smooth function
$\alpha  \in   C^\infty   (V, (0, \infty  ))$:
\begin{enumerate}
 \item $\Gamma$, as a periodic orbit of the perturbed dynamical system
$$\dfrac{du}{dt}= \nu (u) (\nabla H(u) \times \nabla C(u)) + \alpha (u)(C(u) - c) [\nabla H(u) \times (\nabla H(u) \times \nabla C(u))] ,$$
$u \in   V $, is orbitally phase asymptotically stable, with respect to perturbations
in $V$, along the invariant manifold $H^{-1} (\{h\})$.

\item $\Gamma$, as a periodic orbit of the perturbed dynamical system
$$\dfrac{du}{dt}= \nu (u) (\nabla H(u) \times \nabla C(u)) - \alpha (u)(C(u) - c) [\nabla H(u) \times (\nabla H(u) \times \nabla C(u))] ,$$
$u \in   V $, is unstable.
\end{enumerate}
\end{enumerate}
\end{theorem}

\section{Asymptotic stabilization of periodic orbits of the Rabinovich system}
In order to stabilize globally asymptotically an arbitrary periodic orbit from a family of periodic orbits of a generic three-dimensional dynamical system, we consider another class of perturbations.

The dynamics is located on the common level set of the Casimir and the Hamiltonian of the Poisson manifold and it is happening in any tubular neighborhood of an arbitrary periodic orbit of the 
studied system too. In each of the two cases above the dynamics will remain on the surface Hamiltonian=constant, respectively Casimir=constant and we can obtain asymptotic stabilization only on 
the invariant surface where the periodic orbit that we want to stabilize asymptotically is located. If we remove this constrains, the perturbed system will have not any constant of motion and we 
can obtain the global asymptotic stabilization in an entire tubular neighborhood around the periodic orbit considered. In this case remain dynamically invariant only the periodic orbit and the 
two common level surfaces whose intersection represents this periodic orbit considered.
 In order to do this, we do not need any parameterization  of the periodic orbit of the studied system and the same perturbation stabilize 
 globally asymptotically each conexis components if there are more than one for the considered periodic orbit in a 
 tubular neighborhood around it.

\begin{remark}\label{rem}
Let
\begin{equation}\label{eq}
 \dfrac{du}{dt}= \nu (u) (\nabla H(u) \times \nabla C(u)) ,\ u \in   U
 \end{equation}
be a three-dimensional dynamical system defined on an open subset $U \subseteq\R^3$ , with
$H, C, \nu  \in   C^\infty   (U, \R)$ given smooth real functions, such that $H$ and $C$ are functionally
independent on an open subset $V \subseteq U$.

Suppose there exists $\Gamma  \subset V$ a periodic orbit of \eqref{eq}. If $\Gamma  \subseteq  (H, C)^{-1} (\{(h, c)\})$,
where $(h, c) \in   \R^2$ is a regular value for the map $(H, C) : U \to  \R^2$, then if $c$ and $h$ are regular values of the maps 
$C : U \to  \R$, respectively $H : U \to\R$, then for every smooth functions
$\alpha,\beta  \in   C^\infty   (V, (0, \infty  ))$ we have that $\Gamma$ , as a periodic orbit of the perturbed dynamical system
\begin{align*}
\dfrac{du}{dt}=& \nu (u) (\nabla H(u) \times \nabla C(u)) \\
&- \alpha (u)(H(u) - h) [\nabla C(u) \times (\nabla H(u) \times \nabla C(u))] \\
&+\beta (u)(C(u) - c) [\nabla H(u) \times (\nabla H(u) \times \nabla C(u))] ,
\end{align*}
$u \in   V $, is orbitally phase asymptotically stable, with respect to perturbations in $V$.
\end{remark}

Now, for the case of the Rabinovich system, we can apply the Theorem \ref{thm} and the Remark \ref{rem} to obtain three classes of perturbations, namely, one
which keeps dynamically invariant the Hamiltonian $H_\beta$,  one which keeps dynamically invariant the Casimir function $C$, and respectively one which keeps not dynamically invariant 
 the Hamiltonian $H_\beta$ nor the Casimir $C$. 
In all cases, the periodic orbit of the Rabinovich system will remain a periodic orbit for the perturbed system too.

Let us observe that the maximal set where $\nabla H_\beta$ and $\nabla C$ are linearly independent, is the open
set given by the complement of the set of equilibrium points of \eqref{sys}, namely
\begin{equation}
V :=\R^3 \setminus \{\{(x, 0, \beta) : x \in \R \}\cup \{(0, y, −\beta) : y \in \R\}\cup \{(0, 0, z) : z \in \R \}\}
\end{equation}
because $\nabla H_\beta(u)$, $\nabla C(u)$ and $X(u)$ are linearly independent if and only if $u$ is not an equilibrium point of the dynamical system.

Recall from \cite{razvan} that there exists an open and dense subset $S$ of the image of the map
$(H_\beta, C) : \R^3 \to\R^2$ , such that each fiber of any element $(h, c)$ from $S$, corresponds
to periodic orbits of Rabinovich’s equations. Moreover, any such element is a regular value of
$(H_\beta, C)$, as well as its components for the corresponding maps, $H_\beta$ and respectively $C$.

In \cite{razvan} we find the subsets from semialgebraic splitting of all subsets of the image of the energy-Casimir map $Im(H_{\beta},C)$. 
In particular, the subsets of $S$ described in terms of the image of equilibria of the Rabinovich system
through the map $(H_{\beta},C)$ that provide periodic orbits, are:
\begin{enumerate}[1.]
\item $\beta\ne0$

\begin{enumerate}
\item  $c>\beta^2$
\begin{enumerate}[i)]
\item $\Sigma_{(3,\beta)\leftrightarrow (3,\beta)}^{(s,-)\leftrightarrow (s,+)}=\left\{(h,c)\in\R^2:\beta^2c>h^2;c>\beta^2\right\},$
\item $\Sigma_{(3,\beta)}^{(s,-)}=\left\{(h,c)\in\R^2:\beta^2c=h^2;c>\beta^2\right\},$
\item $\Sigma_{(3,\beta)\rightarrow (2,\beta)}^{(s,\pm)\rightarrow (u,*)}=\left\{(h,c)\in\R^2:2h-\beta^2<c<\dfrac{h^2}{\beta^2}; h>\beta^2\right\},$
\item $\Sigma_{(2,\beta)\rightarrow}^{(u,*)\rightarrow}=\left\{(h,c)\in\R^2:\beta^2<c<2h-\beta^2\right\}.$
\end{enumerate}
\item $c=\beta^2$
\begin{enumerate}[i)]
\item $\Sigma_{(1,\beta)\rightarrow (2,\beta)}^{u\rightarrow (u,0)}=\left\{(h,c)\in\R^2:-\beta^2<h<\beta^2;c=\beta^2\right\},$
\item $\Sigma_{(2,\beta)\rightarrow}^{(u,0)\rightarrow}=\left\{(h,c)\in\R^2:h>\beta^2;c=\beta^2\right\}.$
\end{enumerate}

\item $0<c<\beta^2$
\begin{enumerate}[i)]
\item $\Sigma_{(1,\beta)\rightarrow (3,\beta)}^{(s,-)\rightarrow (u,-)}=\left\{(h,c)\in\R^2:-\beta^2<h<0;\max \{-2h-\beta^2,0\}<c<\dfrac{h^2}{\beta^2}\right\},$
\item $\Sigma_{(3,\beta)\rightarrow (3,\beta)}^{(u,-)\rightarrow (u,+)}=\left\{(h,c)\in\R^2:-\beta^2<h<\beta^2;\dfrac{h^2}{\beta^2}<c<\beta^2\right\}$
\item $\Sigma_{(3,\beta)\rightarrow}^{(u,+)\rightarrow}=\left\{(h,c)\in\R^2:h>0;0<c<\min\{\beta^2,\dfrac{h^2}{\beta^2}\}\right\}.$
\end{enumerate}

\item $c=0$
\begin{enumerate}[i)]
\item $\Sigma_{(1,\beta)\rightarrow (3,\beta)}^{(s,0)\rightarrow (u,0)}=\left\{(h,c)\in\R^2:-\dfrac{\beta^2}{2}<h<0;c=0\right\},$
\item $\Sigma_{(3,\beta)\rightarrow}^{(u,0)\rightarrow}=\left\{(h,c)\in\R^2:h>0,c=0\right\}.$
\end{enumerate}

\item $c<0$
\begin{enumerate}[ ]
\item $\Sigma_{(1,\beta)\rightarrow}^{(s,+)\rightarrow}=\left\{(h,c)\in\R^2:-2h-\beta^2<c<0\right\}.$
\end{enumerate}
\end{enumerate}
\item $\beta=0$
\begin{enumerate}
\item $c>0$
\begin{enumerate}[i)]
\item $\Sigma_{(3,0)\rightarrow (2,0)}^{(s,*)\rightarrow u}=\left\{(h,c)\in\R^2:h>0;c>2h\right\},$
\item $\Sigma_{(2,0)\rightarrow}^{u\rightarrow}=\left\{(h,c)\in\R^2:0<c<2h\right\}.$
\end{enumerate}

\item $c=0$
\begin{enumerate}[ ]
\item $\Sigma_{(1,0)\rightarrow}^{(s,0)\rightarrow}=\left\{(h,c)\in\R^2:h>0;c=0\right\}.$
\end{enumerate}

\item $c<0$
\begin{enumerate}[ ]
\item $\Sigma_{(1,0)\rightarrow}^{(s,*)\rightarrow}=\left\{(h,c)\in\R^2:-2h<c<0\right\}.$
\end{enumerate}
\end{enumerate}
\end{enumerate}
For more details see \cite{razvan}.

Let $(h, c) \in S$  and let 
\begin{equation}\label{periodic orbits}
\Gamma  \subseteq  (H_\beta, C)^{-1} (\{(h, c)\})
\end{equation}
be a
periodic orbit of the dynamical system \eqref{sys1}. Then by Theorem \ref{thm}, the following conclusions hold true.

\begin{theorem}

\begin{enumerate}
 \item For every smooth function $a\in C^\infty\left( V,(0,\infty)\right)$:

\begin{enumerate}
 \item $\Gamma$, as a periodic orbit of the dynamical system

$$\left\{
\begin{array}{l}
 \dot x=yz+\beta  y-a(x,y,z) \left(\dfrac{x^2+y^2}{4}-\beta  z-h\right) \left(x y^2+2 x z^2-2 \beta  x z\right) \\
\dot y= -xz+\beta  x-a(x,y,z) \left(\dfrac{x^2+y^2}{4}-\beta  z-h\right) \left(y x^2+2 y z^2+2 \beta  y z\right) \\
 \dot z=x y-a(x,y,z) \left(\dfrac{x^2+y^2}{4} -\beta  z-h\right) \left(-\beta  x^2+z x^2-\beta  y^2-y^2 z\right) \\
\end{array}
\right.$$
$u\in V$, is orbitally phase asymptotically stable, with respect to perturbations in $V$, along the invariant manifold
$$C^{-1}(\{c\})=\left\{(x,y,z)\in{\R}^3\left|\dfrac{1}{2}(-x^2+y^2)+z^2=c\right.\right\}$$

\item $\Gamma$, as a periodic orbit of the dynamical system

$$\left\{
\begin{array}{l}
 \dot x=yz+\beta  y+a(x,y,z) \left(\dfrac{x^2+y^2}{4}-\beta  z-h\right) \left(x y^2+2 x z^2-2 \beta  x z\right) \\
\dot y= -xz+\beta  x+a(x,y,z) \left(\dfrac{x^2+y^2}{4}-\beta  z-h\right) \left(y x^2+2 y z^2+2 \beta  y z\right) \\
 \dot z=x y+a(x,y,z) \left(\dfrac{x^2+y^2}{4} -\beta  z-h\right) \left(-\beta  x^2+z x^2-\beta  y^2-y^2 z\right) \\
\end{array}
\right.$$
$u\in V$, is unstable.
\end{enumerate}

\item For every smooth function $a\in C^\infty\left( V,(0,\infty)\right)$:

\begin{enumerate}
 \item $\Gamma$, as a periodic orbit of the dynamical system

$$\left\{
\begin{array}{l}
 \dot x=yz+\beta  y+a(x,y,z)  \left(z^2+\dfrac{y^2-x^2}{2} -c\right) \left(x \beta ^2-x z \beta +\dfrac{x y^2}{2}\right)\\
\dot y= -xz+\beta  x+a(x,y,z) \left(z^2+\dfrac{y^2-x^2}{2}-c \right)\left(-y z \beta -\beta ^2 y-\dfrac{x^2 y}{2}\right)  \\
\dot z= x y+a(x,y,z)\left(z^2+\dfrac{y^2-x^2}{2}-c\right)\dfrac{\beta  x^2-z x^2-\beta  y^2-y^2 z}{2}  \\
\end{array}
\right.$$
$u\in V$, is orbitally phase asymptotically stable, with respect to perturbations in $V$, along the invariant manifold
$$H_{\beta}^{-1}(\{h\})=\left\{(x,y,z)\in{\R}^3\left|\dfrac{1}{4}(x^2+y^2)-\beta z=h\right.\right\}$$

 \item $\Gamma$, as a periodic orbit of the dynamical system
 $$\left\{
\begin{array}{l}
 \dot x=yz+\beta  y-a(x,y,z)  \left(z^2+\dfrac{y^2-x^2}{2} -c\right) \left(x \beta ^2-x z \beta +\dfrac{x y^2}{2}\right)\\
\dot y= -xz+\beta  x-a(x,y,z) \left(z^2+\dfrac{y^2-x^2}{2}-c \right)\left(-y z \beta -\beta ^2 y-\dfrac{x^2 y}{2}\right)  \\
\dot z= x y-a(x,y,z)\left(z^2+\dfrac{y^2-x^2}{2}-c\right)\dfrac{\beta  x^2-z x^2-\beta  y^2-y^2 z}{2}  \\
\end{array}
\right.$$
$u\in V$, is unstable.
 \end{enumerate}

 \item For every smooth functions $a,\ b\in C^\infty\left( V,(0,\infty)\right)$:
$$\left\{
\begin{split}
\dot x=& yz+\beta  y - a(x,y,z)\left(\frac{x^2+y^2}{4} -\beta  z-h\right) \left(x y^2+2 x z^2-2 \beta  x z\right)  \\
       &+ b(x,y,z) \left(z^2+\frac{y^2-x^2}{2}-c \right)\left(x \beta ^2-x z \beta +\frac{x y^2}{2}\right) \\
\dot y=& -xz+\beta  x-a(x,y,z)\left(\frac{x^2+y^2}{4}-\beta  z-h\right) \left(y x^2+2 y z^2+2 \beta  y z\right)  \\
       &+ b(x,y,z)\left(z^2+\frac{y^2-x^2}{2}-c\right)\left(-y z \beta -\beta ^2 y-\frac{x^2 y}{2}\right)\\
\dot z=&x y-a(x,y,z)\left(\frac{x^2+y^2}{4} -\beta  z-h\right) \left(-\beta  x^2+z x^2-\beta  y^2-y^2 z\right)  \\
       &+ b(x,y,z)\left(z^2+\frac{y^2-x^2}{2} -c\right)\frac{\beta  x^2-z x^2-\beta  y^2-y^2 z}{2}
\end{split}
\right.
$$
 $u\in V$, is orbitally phase asymptotically stable.
\end{enumerate}
\end{theorem}

In the following tables we put toghether in all cases mentioned above the perodic orbit of the perturbated system and the perodic orbit of the Rabinovich system. 
This way we have the complete picture of how we can asymptotically stabilize the periodic orbits of the Rabinovich system.

\begin{landscape}
\begin{figure}[H]
\begin{tabular}{|c|c|c|}
\hline
$\Gamma$ is asymptotically& $\Gamma$ is asymptotically &  $\Gamma$ is asymptotically\\
stable along $C^{-1}(\{c\})$ & stable along $C^{-1}(\{c\})$ and $C^{-1}(\{c'\})$ &  unstable along $C^{-1}(\{c\})$\\
\hline
\includegraphics*{ac1ai.%
eps}&
\includegraphics*{aac1ai.%
eps}
&
\includegraphics*{ic1ai.%
eps}
 \\
 \hline
$\Gamma$ is asymptotically& $\Gamma$ is asymptotically &  $\Gamma$ is orbitally phase\\
stable along $H_{\beta}^{-1}(\{h\})$ & unstable along $H_{\beta}^{-1}(\{h\})$ &  asymptotically stable\\
\hline

\includegraphics*{ah1ai.%
eps}&

\includegraphics*{ih1ai.%
eps}
&
\includegraphics*{tc1ai.%
eps}
 \\
 \hline
\end{tabular}
\centering\caption{Asymptotic stabilization of periodic orbits \eqref{periodic orbits} in the case $\beta\ne0,c>\beta^2;\ c>\frac{h^2}{\beta^2}$}
\label{fig1}
\end{figure}\label{fig1}
\end{landscape}

\begin{landscape}
\begin{figure}[H]
\begin{tabular}{|c|c|c|}
\hline
$\Gamma$ is asymptotically& $\Gamma$ is asymptotically &  $\Gamma$ is asymptotically\\
stable along $C^{-1}(\{c\})$ & stable along $C^{-1}(\{c\})$ and $C^{-1}(\{c'\})$ &  unstable along $C^{-1}(\{c\})$\\
\hline
\includegraphics*{ac1aii.%
eps}&
\includegraphics*{aac1aii.%
eps}
&
\includegraphics*{ic1aii.%
eps}
 \\
\hline
$\Gamma$ is asymptotically& $\Gamma$ is asymptotically &  $\Gamma$ is orbitally phase\\
stable along $H_{\beta}^{-1}(\{h\})$ & unstable along $H_{\beta}^{-1}(\{h\})$ &  asymptotically stable\\
\hline

\includegraphics*{ah1aii.%
eps}&

\includegraphics*{ih1aii.%
eps}
&
\includegraphics*{tc1aii.%
eps}
 \\
 \hline
\end{tabular}
\centering\caption{Asymptotic stabilization of periodic orbits \eqref{periodic orbits} in the case $\beta\ne0,\beta^2<c=\frac{h^2}{\beta^2}$}
\label{fig2}
\end{figure}\label{fig2}
\end{landscape}

\begin{landscape}
\begin{figure}[H]
\begin{tabular}{|c|c|c|}
\hline
$\Gamma$ is asymptotically& $\Gamma$ is asymptotically &  $\Gamma$ is asymptotically\\
stable along $C^{-1}(\{c\})$ & stable along $C^{-1}(\{c\})$ and $C^{-1}(\{c'\})$ &  unstable along $C^{-1}(\{c\})$\\
\hline
\includegraphics*{ac1aiii.%
eps}&
\includegraphics*{aac1aiii.%
eps}
&
\includegraphics*{ic1aiii.%
eps}
 \\
 \hline
$\Gamma$ is asymptotically& $\Gamma$ is asymptotically &  $\Gamma$ is orbitally phase\\
stable along $H_{\beta}^{-1}(\{h\})$ & unstable along $H_{\beta}^{-1}(\{h\})$ &  asymptotically stable\\
\hline

\includegraphics*{ah1aiii.%
eps}&

\includegraphics*{ih1aiii.%
eps}
&
\includegraphics*{tc1aiii.%
eps}
 \\
 \hline
\end{tabular}
\centering\caption{Asymptotic stabilization of periodic orbits \eqref{periodic orbits} in the case $\beta\ne0,2h-\beta^2<c<\frac{h^2}{\beta^2},h>\beta^2$}
\label{fig3}
\end{figure}\label{fig3}
\end{landscape}

\begin{landscape}
\begin{figure}[H]
\begin{tabular}{|c|c|c|}
\hline
$\Gamma$ is asymptotically& $\Gamma$ is asymptotically &  $\Gamma$ is asymptotically\\
stable along $C^{-1}(\{c\})$ & stable along $C^{-1}(\{c\})$ and $C^{-1}(\{c'\})$ &  unstable along $C^{-1}(\{c\})$\\
\hline
\includegraphics*{ac1aiv.%
eps}&
\includegraphics*{aac1aiv.%
eps}
&
\includegraphics*{ic1aiv.%
eps}
 \\
 \hline
$\Gamma$ is asymptotically& $\Gamma$ is asymptotically &  $\Gamma$ is orbitally phase\\
stable along $H_{\beta}^{-1}(\{h\})$ & unstable along $H_{\beta}^{-1}(\{h\})$ &  asymptotically stable\\
\hline

\includegraphics*{ah1aiv.%
eps}&

\includegraphics*{ih1aiv.%
eps}
&
\includegraphics*{tc1aiv.%
eps}
 \\
 \hline
\end{tabular}
\centering\caption{Asymptotic stabilization of periodic orbits \eqref{periodic orbits} in the case $\beta\ne0,\beta^2<c<2h-\beta^2$}
\label{fig4}
\end{figure}\label{fig4}
\end{landscape}

\begin{landscape}
\begin{figure}[H]
\begin{tabular}{|c|c|c|}
\hline
\multicolumn{2}{|c|}{$\Gamma$ is asymptotically table along $C^{-1}(\{c\})$} &  $\Gamma$ is asymptotically unstable along $C^{-1}(\{c\})$\\
\hline
\multicolumn{2}{|c|}{\includegraphics*{ac1bi.%
eps}}
&
\includegraphics*{ic1bi.%
eps}
 \\
 \hline
$\Gamma$ is asymptotically& $\Gamma$ is asymptotically &  $\Gamma$ is orbitally phase\\
stable along $H_{\beta}^{-1}(\{h\})$ & unstable along $H_{\beta}^{-1}(\{h\})$ &  asymptotically stable\\
\hline

\includegraphics*{ah1bi.%
eps}&

\includegraphics*{ih1bi.%
eps}
&
\includegraphics*{tc1bi.%
eps}
 \\
 \hline
\end{tabular}
\centering\caption{Asymptotic stabilization of periodic orbits \eqref{periodic orbits} in the case $\beta\ne0,-\beta^2<h<\beta^2,\ c=\beta^2$}
\label{fig5}
\end{figure}\label{fig5}
\end{landscape}

\begin{landscape}
\begin{figure}[H]
\begin{tabular}{|c|c|c|}
\hline
\multicolumn{2}{|c|}{$\Gamma$ is asymptotically stable along $C^{-1}(\{c\})$} &  $\Gamma$ is asymptotically unstable along $C^{-1}(\{c\})$\\
\hline
\multicolumn{2}{|c|}{\includegraphics*{ac1bii.%
eps}}
&
\includegraphics*{ic1bii.%
eps}
 \\
 \hline
$\Gamma$ is asymptotically& $\Gamma$ is asymptotically &  $\Gamma$ is orbitally phase\\
stable along $H_{\beta}^{-1}(\{h\})$ & unstable along $H_{\beta}^{-1}(\{h\})$ &  asymptotically stable\\
\hline

\includegraphics*{ah1bii.%
eps}&

\includegraphics*{ih1bii.%
eps}
&
\includegraphics*{tc1bii.%
eps}
 \\
 \hline
\end{tabular}
\centering\caption{Asymptotic stabilization of periodic orbits \eqref{periodic orbits} in the case $\beta\ne0,h>\beta^2,c=\beta^2$}
\label{fig6}
\end{figure}\label{fig6}
\end{landscape}

\begin{landscape}
\begin{figure}[H]
\begin{tabular}{|c|c|c|}
\hline
$\Gamma$ is asymptotically& $\Gamma$ is asymptotically &  $\Gamma$ is asymptotically\\
stable along $C^{-1}(\{c\})$ & stable along $C^{-1}(\{c\})$ and $C^{-1}(\{c'\})$ &  unstable along $C^{-1}(\{c\})$\\
\hline
\includegraphics*{ac1ci.%
eps}&
\includegraphics*{aac1ci.%
eps}
&
\includegraphics*{ic1ci.%
eps}
 \\
 \hline
$\Gamma$ is asymptotically& $\Gamma$ is asymptotically &  $\Gamma$ is orbitally phase\\
stable along $H_{\beta}^{-1}(\{h\})$ & unstable along $H_{\beta}^{-1}(\{h\})$ &  asymptotically stable\\
\hline

\includegraphics*{ah1ci.%
eps}&

\includegraphics*{ih1ci.%
eps}
&
\includegraphics*{tc1ci.%
eps}
 \\
 \hline
\end{tabular}
\centering\caption{Asymptotic stabilization of periodic orbits \eqref{periodic orbits} in the case $\beta\ne0,\max(-2h-\beta^2,0)<c<\frac{h^2}{\beta^2},\ -\beta^2<h<0$}
\label{fig7}
\end{figure}\label{fig7}
\end{landscape}

\begin{landscape}
\begin{figure}[H]
\begin{tabular}{|c|c|c|}
\hline
$\Gamma$ is asymptotically& $\Gamma$ is asymptotically &  $\Gamma$ is asymptotically\\
stable along $C^{-1}(\{c\})$ & stable along $C^{-1}(\{c\})$ and $C^{-1}(\{c'\})$ &  unstable along $C^{-1}(\{c\})$\\
\hline
\includegraphics*{ac1cii.%
eps}&
\includegraphics*{aac1cii.%
eps}
&
\includegraphics*{ic1cii.%
eps}
 \\
 \hline
$\Gamma$ is asymptotically& $\Gamma$ is asymptotically &  $\Gamma$ is orbitally phase\\
stable along $H_{\beta}^{-1}(\{h\})$ & unstable along $H_{\beta}^{-1}(\{h\})$ &  asymptotically stable\\
\hline
\includegraphics*{ah1cii.%
eps}&
\includegraphics*{ih1cii.%
eps}
&
\includegraphics*{tc1cii.%
eps}
 \\
 \hline
\end{tabular}
\centering\caption{Asymptotic stabilization of periodic orbits \eqref{periodic orbits} in the case $\beta\ne0,-\beta^2<h<\beta^2,\ \frac{h^2}{\beta^2}<c<\beta^2$}
\label{fig8}
\end{figure}\label{fig8}
\end{landscape}

\begin{landscape}
\begin{figure}[H]
\begin{tabular}{|c|c|c|}
\hline
$\Gamma$ is asymptotically& $\Gamma$ is asymptotically &  $\Gamma$ is asymptotically\\
stable along $C^{-1}(\{c\})$ & stable along $C^{-1}(\{c\})$ and $C^{-1}(\{c'\})$ &  unstable along $C^{-1}(\{c\})$\\
\hline
\includegraphics*{ac1ciii.%
eps}&
\includegraphics*{aac1ciii.%
eps}
&
\includegraphics*{ic1ciii.%
eps}
 \\
 \hline
$\Gamma$ is asymptotically& $\Gamma$ is asymptotically &  $\Gamma$ is orbitally phase\\
stable along $H_{\beta}^{-1}(\{h\})$ & unstable along $H_{\beta}^{-1}(\{h\})$ &  asymptotically stable\\
\hline
\includegraphics*{ah1ciii.%
eps}
&
\includegraphics*{ih1ciii.%
eps}
&
\includegraphics*{tc1ciii.%
eps}
 \\
 \hline
\end{tabular}
\centering\caption{Asymptotic stabilization of periodic orbits \eqref{periodic orbits} in the case $\beta\ne0, 0<c<\min(\beta^2,\frac{h^2}{\beta^2}),h>0$}
\label{fig9}
\end{figure}\label{fig9}
\end{landscape}

\begin{landscape}
\begin{figure}[H]
\begin{tabular}{|c|c|c|}
\hline
\multicolumn{2}{|c|}{$\Gamma$ is asymptotically stable along $C^{-1}(\{c\})$} &  $\Gamma$ is asymptotically unstable along $C^{-1}(\{c\})$\\
\hline
\multicolumn{2}{|c|}{\includegraphics*{ac1di.%
eps}}
&
\includegraphics*{ic1di.%
eps}
 \\
 \hline
$\Gamma$ is asymptotically& $\Gamma$ is asymptotically &  $\Gamma$ is orbitally phase\\
stable along $H_{\beta}^{-1}(\{h\})$ & unstable along $H_{\beta}^{-1}(\{h\})$ &  asymptotically stable\\
\hline
\includegraphics*{ah1di.%
eps}
&
\includegraphics*{ih1di.%
eps}
&
\includegraphics*{tc1di.%
eps}
 \\
 \hline
\end{tabular}
\centering\caption{Asymptotic stabilization of periodic orbits \eqref{periodic orbits} in the case $\beta\ne0,-\frac{\beta^2}{2}<h<0,\ c=0$}
\label{fig10}
\end{figure}\label{fig10}
\end{landscape}

\begin{landscape}
\begin{figure}[H]
\begin{tabular}{|c|c|c|}
\hline
\multicolumn{2}{|c|}{$\Gamma$ is asymptotically stable along $C^{-1}(\{c\})$} &  $\Gamma$ is asymptotically  unstable along $C^{-1}(\{c\})$\\
\hline
\multicolumn{2}{|c|}{\includegraphics*{ac1dii.%
eps}}
&
\includegraphics*{ic1dii.%
eps}
 \\
 \hline
$\Gamma$ is asymptotically& $\Gamma$ is asymptotically &  $\Gamma$ is orbitally phase\\
stable along $H_{\beta}^{-1}(\{h\})$ & unstable along $H_{\beta}^{-1}(\{h\})$ &  asymptotically stable\\
\hline
\includegraphics*{ah1dii.%
eps}
&
\includegraphics*{ih1dii.%
eps}
&
\includegraphics*{tc1dii.%
eps}
 \\
 \hline
\end{tabular}
\centering\caption{Asymptotic stabilization of periodic orbits \eqref{periodic orbits} in the case $\beta\ne0,h>0,\ c=0$}
\label{fig11}
\end{figure}\label{fig11}
\end{landscape}

\begin{landscape}
\begin{figure}[H]
\begin{tabular}{|c|c|c|}
\hline
$\Gamma$ is asymptotically& $\Gamma$ is asymptotically &  $\Gamma$ is asymptotically\\
stable along $C^{-1}(\{c\})$ & stable along $C^{-1}(\{c\})$ and $C^{-1}(\{c'\})$ &  unstable along $C^{-1}(\{c\})$\\
\hline
\includegraphics*{ac1e.%
eps}
&
\includegraphics*{aac1e.%
eps}
&
\includegraphics*{ic1e.%
eps}
 \\
 \hline
$\Gamma$ is asymptotically& $\Gamma$ is asymptotically &  $\Gamma$ is orbitally phase\\
stable along $H_{\beta}^{-1}(\{h\})$ & unstable along $H_{\beta}^{-1}(\{h\})$ &  asymptotically stable\\
\hline
\includegraphics*{ah1e.%
eps}
&
\includegraphics*{ih1e.%
eps}
&
\includegraphics*{tc1e.%
eps}
 \\
 \hline
\end{tabular}
\centering\caption{Asymptotic stabilization of periodic orbits \eqref{periodic orbits} in the case $\beta\ne0,-2h-{\beta^2}<c<0$}
\label{fig12}
\end{figure}\label{fig12}
\end{landscape}

\begin{landscape}
\begin{figure}[H]
\begin{tabular}{|c|c|c|}
\hline
$\Gamma$ is asymptotically& $\Gamma$ is asymptotically &  $\Gamma$ is asymptotically\\
stable along $C^{-1}(\{c\})$ & stable along $C^{-1}(\{c\})$ and $C^{-1}(\{c'\})$ &  unstable along $C^{-1}(\{c\})$\\
\hline
\includegraphics*{ac0ai.%
eps}
&
\includegraphics*{aac0ai.%
eps}
&
\includegraphics*{ic0ai.%
eps}
 \\
 \hline
$\Gamma$ is asymptotically& $\Gamma$ is asymptotically &  $\Gamma$ is orbitally phase\\
stable along $H_{\beta}^{-1}(\{h\})$ & unstable along $H_{\beta}^{-1}(\{h\})$ &  asymptotically stable\\
\hline
\includegraphics*{ah0ai.%
eps}
&
\includegraphics*{ih0ai.%
eps}
&
\includegraphics*{tc0ai.%
eps}
 \\
 \hline
\end{tabular}
\centering\caption{Asymptotic stabilization of periodic orbits \eqref{periodic orbits} in the case $\beta=0,h>0,c>2h$}
\label{fig13}
\end{figure}\label{fig13}
\end{landscape}

\begin{landscape}
\begin{figure}[H]
\begin{tabular}{|c|c|c|}
\hline
$\Gamma$ is asymptotically& $\Gamma$ is asymptotically &  $\Gamma$ is asymptotically\\
stable along $C^{-1}(\{c\})$ & stable along $C^{-1}(\{c\})$ and $C^{-1}(\{c'\})$ &  unstable along $C^{-1}(\{c\})$\\
\hline
\includegraphics*{ac0aii.%
eps}
&
\includegraphics*{aac0aii.%
eps}
&
\includegraphics*{ic0aii.%
eps}
 \\
 \hline
$\Gamma$ is asymptotically& $\Gamma$ is asymptotically &  $\Gamma$ is orbitally phase\\
stable along $H_{\beta}^{-1}(\{h\})$ & unstable along $H_{\beta}^{-1}(\{h\})$ &  asymptotically stable\\
\hline
\includegraphics*{ah0aii.%
eps}
&
\includegraphics*{ih0aii.%
eps}
&
\includegraphics*{tc0aii.%
eps}
 \\
 \hline
\end{tabular}
\centering\caption{Asymptotic stabilization of periodic orbits \eqref{periodic orbits} in the case $\beta=0,0<c<2h$}
\label{fig14}
\end{figure}\label{fig14}
\end{landscape}

\begin{landscape}
\begin{figure}[H]
\begin{tabular}{|c|c|c|}
\hline
\multicolumn{2}{|c|}{$\Gamma$ is asymptotically  stable along $C^{-1}(\{c\})$} &  $\Gamma$ is asymptotically  unstable along $C^{-1}(\{c\})$\\
\hline
\multicolumn{2}{|c|}{\includegraphics*{ac0b.%
eps}}
&
\includegraphics*{ic0b.%
eps}
 \\
 \hline
$\Gamma$ is asymptotically& $\Gamma$ is asymptotically & $\Gamma$ is orbitally phase\\
stable along $H_{0}^{-1}(\{h\})$ &unstable along $H_{0}^{-1}(\{h\})$ &  asymptotically stable\\
\hline
\includegraphics*{ah0b.%
eps}
&
\includegraphics*{ih0b.%
eps}
&
\includegraphics*{tc0b.%
eps}
 \\
 \hline
\end{tabular}
\centering\caption{Asymptotic stabilization of periodic orbits \eqref{periodic orbits} in the case $\beta=0,h>0,c=0$}
\label{fig15}
\end{figure}\label{fig15}
\end{landscape}

\begin{landscape}
\begin{figure}[H]
\begin{tabular}{|c|c|c|}
\hline
$\Gamma$ is asymptotically& $\Gamma$ is asymptotically &  $\Gamma$ is asymptotically\\
stable along $C^{-1}(\{c\})$ & stable along $C^{-1}(\{c\})$ and $C^{-1}(\{c'\})$ &  unstable along $C^{-1}(\{c\})$\\
\hline
\includegraphics*{ac0c.%
eps}
&
\includegraphics*{aac0c.%
eps}
&
\includegraphics*{ic0c.%
eps}
 \\
 \hline
$\Gamma$ is asymptotically& $\Gamma$ is asymptotically &  $\Gamma$ is orbitally phase\\
stable along $H_{0}^{-1}(\{h\})$ & unstable along $H_{0}^{-1}(\{h\})$ &  asymptotically stable\\
\hline
\includegraphics*{ah0c.%
eps}
&
\includegraphics*{ih0c.%
eps}
&
\includegraphics*{tc0c.%
eps}
 \\
 \hline
\end{tabular}
\centering\caption{Asymptotic stabilization of periodic orbits \eqref{periodic orbits} in the case $\beta=0,-2h<c<0$}
\label{fig16}
\end{figure}\label{fig16}
\end{landscape}

\section*{References}


\end{document}